\documentclass[11pt, reqno]{amsart}

\usepackage{mathrsfs,bbold}
\usepackage{amsmath, amsthm, amssymb}
\usepackage{enumerate}
\usepackage{dsfont}
\usepackage{color}
\usepackage[a4paper]{geometry}
\usepackage[utf8]{inputenc}
\usepackage{stmaryrd}
\usepackage{titlesec, fancyhdr}
\usepackage{mathabx}
\usepackage{appendix}
\usepackage{todonotes}
\usepackage{hyperref} 

\newcommand{\RR}{\mathbb{R}}

 \newcommand{\mcC}{\mathcal{C}}

\newcommand{\varep}{\varepsilon}
\newcommand{\ep}{\epsilon}

\renewcommand{\leq}{\leqslant}
\renewcommand{\geq}{\geqslant}

\newcommand{\ssk}{\smallskip}

\newcommand{\cdummy}{\cdot}
\newcommand{\mathd}{\mathrm{d}}

\newcommand{\tmop}[1]{\ensuremath{\operatorname{#1}}}

\newcommand{\vertiii}[1]{{\left\vert\kern-0.25ex\left\vert\kern-0.25ex\left\vert #1 
    \right\vert\kern-0.25ex\right\vert\kern-0.25ex\right\vert}}

\newenvironment{Dem}[1][\unskip]{
    \begin{list}
    {\hspace{0.5cm}{\sf \textbf{Proof #1 --}}}
    {   \setlength{\topsep}{0pt}%
        \setlength{\leftmargin}{0pt}%
        \setlength{\rightmargin}{0pt}%
        \setlength{\listparindent}{0pt}%
        \setlength{\itemindent}{0pt}%
        \setlength{\parsep}{0pt}%
        \addtolength{\leftmargin}{20pt}%
        \addtolength{\rightmargin}{0pt}%
	} 
	\item }{\hfill $\rhd$
	\end{list}}

\newtheoremstyle{mystyle}
{3pt}               
{3pt}               
{\it }                      
{}                      
{\sffamily\bfseries}             
{}                      
{0.5em}                 
{\llap{#2. }#1{$\;$ --}}

\theoremstyle{mystyle}

\newtheorem{thm}{Theorem}
\newtheorem*{thm*}{Theorem}

\newtheorem{cor}[thm]{\hspace{-0.15cm}  {Corollary} }
\newtheorem{lem}[thm]{\hspace{-0.2cm}  {Lemma} }
\newtheorem{prop}[thm]{\hspace{-0.2cm} {Proposition}}
\newtheorem*{defn*} {Definition}
\newtheorem*{prop*} {Proposition}
\newtheorem*{lem*} {Lemma}
\newtheorem*{cor*} {Corollary}

\newtheoremstyle{mystyle2}
{3pt}               
{3pt}               
{\it }                      
{}                      
{\sffamily\bfseries}             
{}                      
{0.5em}                 
{\llap{#2 }#1{$\;$ --}}
\theoremstyle{mystyle2}

\newtheorem*{definition*}{Definition}

\titleformat{\section}[block]
{\filcenter\normalfont\sffamily\bfseries\Large}{{\hspace{-0.7cm}}\thesection \vspace{0.3cm}}{0.75em}{}

\titleformat{\subsection}[block]
{\filcenter\normalfont\sffamily\bfseries\Large}{{\hspace{-0.7cm}}\thesubsection \vspace{0.3cm}}{0.75em}{}

\titleformat{\subsubsection}[block]
{\normalfont\sffamily\bfseries}{\hspace{-1.2cm}\thesubsubsection  \vspace{0.3cm}}{.75em}{}
\titlespacing{\subsection}{0pc}{1.5ex plus .1ex minus .2ex}{0pc}

\numberwithin{subsection}{section}
\numberwithin{subsubsection}{subsection}

\let\oldtocsection=\tocsection
\let\oldtocsubsection=\tocsubsection
\let\oldtocsubsubsection=\tocsubsubsection
\renewcommand{\tocsection}[2]{\hspace{0em}\oldtocsection{#1}{#2}}
\renewcommand{\tocsubsection}[2]{\hspace{1em}\oldtocsubsection{#1}{#2}}
\renewcommand{\tocsubsubsection}[2]{\hspace{2em}\oldtocsubsubsection{#1}{#2}}
\begin{document}

\vspace*{3ex minus 1ex}
\begin{center}
\huge\sffamily{Quasilinear generalized parabolic Anderson model equation}
\end{center}
\vskip 5ex minus 1ex

\begin{center}
{\sf I. BAILLEUL\footnote{I. Bailleul thanks the U.B.O. for their hospitality.} and A. DEBUSSCHE\footnote{I. Bailleul and A. Debussche benefit from the support of the french government ``Investissements d'Avenir'' program ANR-11-LABX-0020-01.} and M. HOFMANOV\'A\footnote{M. Hofmanov\'a gratefully acknowledges the financial support by the DFG via Research Unit FOR 2402.}}
\end{center}

\vspace{1cm}

\begin{center}
\begin{minipage}{0.8\textwidth}
\renewcommand\baselinestretch{0.7} \rmfamily {\scriptsize {\noindent \textsc{\textbf{Abstract.}}}  
We present in this note a local in time well-posedness result for the singular $2$-dimensional quasilinear generalized parabolic Anderson model equation
$$
\partial_t u - a(u)\Delta u = g(u)\xi
$$
The key idea of our approach is a simple transformation of the equation which allows to treat the problem as a semilinear problem. The analysis is done within the elementary setting of paracontrolled calculus.}
\end{minipage}
\end{center}

\bigskip
\bigskip

\section{Introduction}
\label{SectionIntro}

Tremendous progress has been made recently in the application of rough path ideas to the construction of solutions to singular stochastic partial differential equations (PDEs) driven by time/space rough perturbations, in particular, using Hairer's theory of regularity structures \cite{hairer_theory_2014} and the tools of paracontrolled calculus introduced by Gubinelli, Imkeller and Perkowski in \cite{GIP}. We refer the reader to the works \cite{ AllezChouk,BailleulBernicot2,BHZ, HP15, HQ,  MW} for a tiny sample of the exponentially growing literature on the subject. The (generalised) parabolic Anderson model equation itself was studied from both points of view in different settings in \cite{ BailleulBernicot1, BailleulBernicot2, BailleulBernicotFrey, GIP, hairer_theory_2014, hairer_simple_2015}. The class of equations studied so far in the literature centers around semilinear problems with nonlinear lower order terms that are not well defined in the classical sense. We investigate in the present paper the possibility of extending these methods towards a quasilinear setting, that is, towards problems with nonlinear dependence on the solution in the leading order term. The first result in this direction was obtained very recently by Otto and Weber in their work \cite{OttoWeber}, in which they study the $(1+1)$-dimensional \textit{time/space periodic} equation
$$
\partial_t u - {\sf P}\big(a (u) \Delta u\big) = {\sf P}\big(g (u) \xi\big),
$$
with a mildy irregular time/space periodic noise $\xi$ of parabolic H\"older regularity $(\alpha-2)$, for $\frac{2}{3}<\alpha<1$ -- ${\sf P}$ stands here for the projection operator on zero spatial mean functions. They develop for that purpose a simplified, parametric, version of regularity structures in the line of Gubinelli's approach to rough differential equations using controlled paths. This approach requires a whole new setting that is described at length in \cite{OttoWeber}. 
It is elegantly rephrased by Furlan and Gubinelli   \cite{GubinelliFurlan} in a work which is independent and simultaneous to the present one. They use a variant of paracontrolled calculus based on paracomposition operators for the study of the \textit{evolution} quasilinear equation 
$$
\partial_t u - a (u) \Delta u = g (u) \xi, \qquad u (0) = u_0, 
$$
where $\xi$ is a space white noise on the $2$-dimensional torus. This is the equation which we study here.

\bigskip

Recall that the zero mean space white noise $\xi$ over the $2$-dimensional torus is almost surely of spatial H\"older regularity $\alpha-2$, for any $\frac{2}{3}<\alpha<1$, and write $(\xi^\epsilon)_{0<\epsilon\leq 1}$ for the family of smoothened noises obtained by convolution of $\xi$ with the heat kernel. See below for the definition of the spatial and parabolic H\"older spaces $C^\alpha$ and $\mcC^\alpha_T$ mentioned in the statement.

\medskip

\begin{thm}   \label{ThmMain}
Let a function $a\in C^3_b$, with values in some compact interval of $(0,\infty)$, and a function $g\in C^3_b$ be given. Let also a regularity exponent $\alpha\in\big(\frac{2}{3},1\big)$ be given, together with an initial condition $u_0\in C^\alpha$. Then there are some diverging constants $c^\epsilon$ and a random time $T$, defined on the same probability space as space white noise, such that the solutions $u^\epsilon$ to the well-posed equations 
\begin{equation*}
\label{EqRegularized}
\partial_t u^\epsilon - a(u^\epsilon)\Delta u^\epsilon = g(u^\epsilon)\,\xi^\epsilon - c^\epsilon\,\left\{\Big(\frac{g'g}{a}\Big)(u^\epsilon) - \Big(\frac{a'g^2}{a^2}\Big)(u^\epsilon)\right\}
\end{equation*}
with initial value $u_0$, converge, as $\epsilon$ decreases to $0$, almost surely in the parabolic H\"older space $\mcC^\alpha_T$ to a limit element $u\in\mcC^\alpha_T$, unique solution of the paracontrolled singular equation
\begin{equation}  
\label{EqQuasigPAM}
\partial_t u - a (u) \Delta u = g (u) \xi, \qquad u (0) = u_0. 
\end{equation}
\end{thm}

\medskip

A solution to a paracontrolled singular equation is more properly a pair $(u,u')$; the above improper formulation is justified in so far as $u'$ will actually be a function of $u$. This statement is the exact analogue of the main result obtained by Furlan and Gubinelli in \cite{GubinelliFurlan}, using their extension of paracontrolled calculus based on paracomposition operators. The present work makes it clear that the basic tools of paracontrolled analysis are sufficient for the analysis of this equation. Note here the slight improvement over \cite{GubinelliFurlan} in the convergence of $u^\epsilon$ to the limit function $u$, that takes place here in the parabolic H\"older space $\mcC^\alpha_T$ rather than just in $C_TC^\alpha$. Note here that our approach works verbatim if one replaces the operator $a(u)\Delta u$ by $a^{ij}(u)\partial^2_{ij}$, for some matrix-valued function $a(u)$ that is symmetric and uniformly elliptic, and for $u$ taking values in some finite dimensional vector space. Adding a term $b^i(u)\partial_i $ in the dynamics would not cause any trouble in the range of regularity $\alpha\in\big(\frac{2}{3},1\big)$ where we are working. Last, note that in the \textit{scalar-valued case}, solving equation \eqref{EqQuasigPAM} is \textit{equivalent} to solving an equation of the form 
$$
\partial_t v - \Delta\big(b(v)\big) = f(v) \xi,\qquad v(0) = v_0,
$$
after setting $v := A(u)$, with $A$ a primitive of $1/a$, with $b$ the inverse of $A$, and $f = (g/a)\circ b$.

\bigskip

The precise setting of paracontrolled calculus that will be used here in the analysis of the singular equation \eqref{EqQuasigPAM} is detailed in Section \ref{SectionSetting}, where the proof of Theorem \ref{ThmMain} is given in three steps. A number of elementary results have been put aside in Appendix.

\vfill
\pagebreak

\noindent \textsf{\textbf{Notations.}} We gather here a number of notations that will be used in the text.   \vspace{0.15cm}

\begin{itemize}
   \item[\textcolor{gray}{$\bullet$}] Let $P$ stand for the heat semigroup associated with the Laplace operator $\Delta$ on the $2$-dimensional torus, $\mathscr{L} := \partial_t - \Delta$ stand for the heat operator, and $\mathscr{L}^{-1}$ stand for the resolution operator of the heat equation $\mathscr{L} u := (\partial_t - \Delta) u  = f$, with null initial condition, given by the formula 
$$
\big(\mathscr{L}^{-1}f\big)(t) := \int_0^t P_{t-s}f_s\,ds,
$$
for a time-dependent distribution $f$. We use a similar notation if $\Delta$ is replaced by another uniformly elliptic operator.   \vspace{0.1cm}

   \item[\textcolor{gray}{$\bullet$}] Given a positive time horizon $T$, a regularity exponent $\alpha$, and  a Banach space $E$, write $C^\alpha_T E$ for $C^\alpha\big([0,T],E\big)$. Given a real regularity exponent $\alpha$, we denote by $C^\alpha$ the spatial H\"older space and by $\mcC^\alpha$ the parabolic H\"older space, both defined for instance in terms of Besov spaces built from the parabolic operator $\mathscr{L}$ -- see e.g. \cite{BailleulBernicot1}. For $\alpha\in (0,2)$, the parabolic space $\mcC^\alpha$, or $\mcC^\alpha_T$, coincides as a set with $C_T C^{\alpha} \cap C_T^{\alpha / 2} L^{\infty}$, and the Besov norm on $\mcC^\alpha$ is equivalent to the elementary norm  
$$
\| \cdummy \|_{\mcC^{\alpha}} := \| \cdummy \|_{C_T C^{\alpha}} + \| \cdummy \|_{C_T^{\alpha / 2} L^{\infty}}.
$$
\end{itemize}

\bigskip

\section{Paracontrolled setting}
\label{SectionSetting}

Let $\xi$ stand for a space white noise on the $2$-dimensional torus. In its simplest form, the multiplication problem raised by an ill-posed product, like the term $g(u)\xi$ in the model $2$-dimensional generalised (PAM) equation 
$$
\mathscr{L}u = g(u)\xi,
$$
is dealt within paracontrolled calculus by looking for solutions $u$ of the equation in an a priori rigidly structured, graded, solution space, whose elements locally look like some reference function built from $\xi$ only by classical means. This approach requires the problem-independent assumption that some product(s) involving only the noise $\xi$ can be  given an analytical sense in some appropriate distribution/function space(s); this is typically done using some probabilistic tools. The datum of $\xi$ and all these distributions defines the \textsf{\textbf{enhanced noise}} $\widehat{\xi}$. In the present setting where the spatial noise $\xi$ is $(\alpha-2)$-H\"older, for $\frac{2}{3}<\alpha<1$, only one extra component needs to be added to $\xi$ to get the enhanced noise, the associated solution space has two levels, and a potential solution $u$ two components $\big(u,u'\big)$. The structure of the elements in the solution space and the datum of the enhanced noise allow for a proper analytical definition fo the product $g(u)\xi$ and show that 
$$
v := Pu_0 + \mathscr{L}^{-1}\big(g(u)\xi\big)
$$
takes values in the solution space; write $\Phi\big(u,u'\big) := \big(v,v'\big)$ for its components, with $\Phi$ a regular fuction of $(u,u'\big)$. The equation 
$$
u = Pu_0 + \mathscr{L}^{-1}\big(g(u)\xi\big),
$$
or rather, 
$$
(u,u'\big) = \Phi\big(u,u'\big)
$$
is then solved on a short time interval using a fixed point argument. The short time horizon is what provides the local contracting character of the map $\Phi$. This scheme works particularly well for an initial condition $u_0\in C^{2\alpha}$, as the term $Pu_0$ can then be inserted in some remainder term; see for instance \cite{BailleulBernicot1}. One needs to adopt a different functional setting for the levels of the solution space to work with an initial condition $u_0$ in $C^\alpha$, see \cite{GIP}. Schauder estimates are used crucially in this reasoning to ensure that $\mathscr{L}^{-1}\big(g(u)\xi\big)$ takes values in the solution space, and one has good quantitative controls on its two levels.

\ssk

The situation gets more complex in the quasilinear setting where the operator 
$$
\mathscr{L}^u := \partial_t - a(u)\Delta
$$
depends itself on the solution $u$ of the equation $\mathscr{L}^uu = g(u)\xi$. Both Otto-Weber \cite{OttoWeber} and Furlan-Gubinelli \cite{GubinelliFurlan} work with parameter dependent operators $\mathscr{L}^b := \partial_t-b\Delta$, for a real-valued positive parameter $b$ ranging in a compact subinterval of $(0,\infty)$, and take profit from $b$-uniform Schauder estimates. The difficulty in their approaches is to get a fixed point reformulation of the equation in an adequate setting. The new rough path-flavoured setting developed at length by Otto and Weber in \cite{OttoWeber} has an elegant counterpart in the relatively short work \cite{GubinelliFurlan} of Furlan and Gubinelli, that requires an extension of paracontrolled calculus via the introduction of paracomposition operators and associated continuity results, mixing seminal works of Alinhac \cite{Alinhac} in the 80's and the basic tools of paracontrolled calculus \cite{GIP}. \textit{We show in the present work that the analysis of the quasilinear generalised} (PAM) \textit{equation \eqref{EqQuasigPAM} can be run efficiently using the elementary paracontrolled calculus, with no need of any new tools.}

\bigskip

The foundations of paracontrolled calculus were laid down in the seminal work \cite{GIP} of Gubinelli, Imkeller and Perkowski, to which we shall refer the reader for a number of facts used here -- see also \cite{BailleulBernicot1, BailleulBernicot2, BailleulBernicotFrey} for extensions. We refer to the book \cite{BCD} of Bahouri, Chemin and Danchin for a gentle introduction to the use of paradifferential calculus in the study of nonlinear PDEs. We shall then freely use the notations $\Pi_fg$ and $\Pi(f,g)$ for the paraproduct of $f$ by $g$ and the corresponding resonant term, defined in terms of Littlewood-Paley decomposition, for any two functions $f$ and $g$ in some spatial H\"older spaces of any regularity exponent. We will denote by $\overline\Pi$ the modified paraproduct on parabolic functions/distributions introduced in \cite{GIP}, formula (36) in Section 5, in which the time fluctuations of the low frequency distribution/function are averaged differently at each space scale. (This modified paraproduct is different from the parabolic paraproduct introduced in \cite{BailleulBernicotFrey}.) The following definition of a paracontrolled distribution will make clear what a ``rigidly structured, graded, solution space'' may look like. \textit{We fix once and for all some regularity exponents $\alpha \in \left(\frac{2}{3}, 1\right)$ and} $\beta\in\left(\frac{2}{3}\vee \frac{\alpha}{2}, \alpha\right)$, and define $X\in C^\alpha$ as the (random) zero spatial mean solution of the equation
$$
-\Delta X = \xi.
$$

\medskip

\begin{definition*}
We define the space ${\sf C}_{\alpha,\beta}(X)$ of \textsf{\textbf{functions paracontrolled by $X$}} as the set of pairs of parabolic functions $ (u, u') \in \mcC^\alpha\times\mcC^\beta$ such that
$$
u^\sharp := u - \overline\Pi_{u'} X \in \mcC^{\alpha}
$$
satisfies
  \begin{equation*}
  \label{eq:expl}
    \sup_{0<t\leq T} t^{\frac{2\beta-\alpha}{2}} \, \big\| u^\sharp\big\|_{C^{2 \beta}} < \infty. 
  \end{equation*}
Setting 
$$
\big\| (u,u') \big\|_{\alpha,\beta} := \| u' \|_{\mcC^\beta} +\| u^{\sharp}\|_{\mcC^\alpha} +  \sup_{0<t\leq T} t^{\frac{2\beta-\alpha}{2}} \, \big\| u^\sharp\big\|_{C^{2\beta}}
$$
turns ${\sf C}_{\alpha,\beta}(X)$ into a Banach space.
\end{definition*}

\medskip

Mention here that different choices can be done for the norm on the space of controlled functions; different purposes may lead to different choices -- see for instance the study of the $2$-dimensional generalised (PAM) equation done in \cite{BailleulBernicot1}. Given a positive time horizon $T$, set 
$$
u_0^T := P_Tu_0,
$$ 
to shorten notations, and recall for future use the bounds
$$
\|u_0^T\|_{C^{2\beta}} \lesssim T^{-\frac{2\beta-\alpha}{2}}\,\|u_0\|_{C^\alpha}
$$
and 
$$
\big\|u_0^T-u_0\big\|_{L^\infty} \lesssim T^\frac{\alpha}{2}\,\|u_0\|_{C^\alpha}.
$$
Our starting point for the analysis of the quasilinear generalised (PAM) equation
$$
\partial_t u -a(u)\Delta u = g(u)\xi
$$
is to rewrite it under the form
\begin{equation}
\label{EqRewritegPAM}
\mathscr{L}^0 u := \partial_t u -a(u_0^T)\Delta u = g(u)\xi + \big(a(u)-a(u_0^T)\big)\Delta u.
\end{equation}
Notice that the term $\big(a(u)-a(u_0^T)\big)\Delta u$ is still part of the leading order operator, so one cannot expect to treat it by perturbation methods. A  suitable paracontrolled ansatz allows however to cancel out the most irregular part of this term and leave us with a remainder that has a better regularity and can then be treated as a perturbation term. This simplification mechanism rests heavily on the fact that since the solution remains ``close'' to its initial value on a small time interval, the difference $a(u)-a(u_0^T)$, and $\big(a(u)-a(u_0^T)\big)\Delta u$ with it, needs to be small in a suitable sense.

\bigskip

\subsection{Fixed point setting}
\label{SubsectionFixedPtSetting}

We rephrase in this section equation \eqref{EqQuasigPAM}, or equivalently equation \eqref{EqRewritegPAM}, as a fixed point problem in ${\sf C}_{\alpha,\beta}(X)$ for a regular map $\Phi$ from ${\sf C}_{\alpha,\beta}(X)$ to itself. This requires that we first make it clear that the a priori ill-defined products $g(u)\xi$ and $\big(a(u)-a(u_0^T)\big)\Delta u$ actually make sense for $u$ paracontrolled by $X$. This holds under the assumption that $\Pi(X,\xi)$ can be properly defined as an element of $C_TC^{2\alpha-2}$; such matters are dealt with in Section \ref{SectionRenormalisation}, from which Theorem \ref{ThmMain} will follow. Given the regularizing properties of the resolution operator of the operator $\mathscr{L}^0$, encoded in the Schauder estimates that we shall use below, terms in $C_TC^{\geq \alpha+ \beta-2}$, where the space $C_TC^{\geq \alpha+ \beta-2}$ is the intersection of $C_TC^{\gamma}$ for 
$\gamma \geq \alpha+ \beta-2$, will be considered as remainders in the analysis of different terms done in this section.

\medskip

The analysis of the term $g(u)\xi$ is done as in Gubinelli, Imkeller and Perkowski's treatement of the generalised (PAM) equation, using paralinearisation of $g(u)$ and the continuity of the correctors 
\begin{equation*}
\begin{split}
& {\sf C}(a,b,\xi) = \Pi\big(\Pi_ab,\xi\big) - a\Pi(b,\xi),   \\
& \overline{\sf C}(a,b,\xi) := \Pi\big(\overline\Pi_ab,\xi\big) - a\Pi(b,\xi),
\end{split}
\end{equation*}
from $C_T C^{\alpha_1}\times C_T C^{\alpha_2}\times C^{\alpha-2}$ to $C_T C^{\alpha_1+\alpha_2 + \alpha-2}$, provided $\alpha_1+\alpha_2 + \alpha > 2$. It gives
\begin{equation}
\label{EqRenormalisation1}
g(u)\xi = \Pi_{g(u)}\xi + \Pi_\xi \big(g(u)\big) + g'(u)u'\,\Pi(X,\xi) + (\star),
\end{equation}
for a remainder $(\star)$ that is the sum of $g'(u)\Pi\big(u^\sharp,\xi\big)$ and an element of $C_T C^{2\alpha + \beta - 2}$, whose norm depends polynomially on $\big\|(u,u')\big\|_{\alpha,\beta}$, and \textit{which depends continuously on} $\xi\in C^{\alpha-2}$. Given the regularity assumption on $u^\sharp$ the term $g'(u)\Pi\big(u^\sharp,\xi\big)$ can only be evaluated in a weighted space.

\medskip

We shall use the following elementary lemma to make clear that the term $\big(a(u)-a(u_0^T)\big)\Delta u$ is well-defined when $u$ is paracontrolled by $X$, and give a description for it up to some remainder term; the proof of the lemma is given in Appendix for completeness. 


\begin{lem}
\label{LemFixedPointSetting}
The following two estimates hold.   \vspace{0.1cm}

\begin{itemize}
   \item[\textcolor{gray}{$\bullet$}] Let $f,g\in C^\beta$, $a\in C^\alpha$ and $b\in C^{\alpha-2}$ be given. Then
   $$
   \Big\| \Pi\big(\Pi_f a,\Pi_g b\big) - fg\,\Pi(a,b)\Big\|_{C^{2\alpha+\beta-2}} \lesssim \|f\|_{C^\beta} \|g\|_{C^\beta} \|a\|_{C^\alpha}\|b\|_{C^{\alpha-2}}.  
   $$   
   
   \item[\textcolor{gray}{$\bullet$}] For $f$ in the parabolic H\"older space $\mcC^\beta$, we have the intertwining continuity estimate
   $$
   \Big\|\mathscr{L}^0\big(\overline\Pi_f X\big) - \Pi_{a(u_0^T)f}\big(-\Delta X\big)\Big\|_{C_TC^{\alpha+\beta-2}} \lesssim \Big(1+T^{-{\frac{2\beta-\alpha}{2}}}\,\|u_0\|_{C^\alpha}\Big)\|f\|_{\mcC^\beta} \|X\|_{C^\alpha}   .
   $$ 
\end{itemize}
\end{lem}

\medskip

We shall use the notation $(\Asterisk)$ for an element of $C_TC^{\geq \alpha+\beta-2}$ which may change from line to line, but \textit{which depends continuously on} $\xi\in C^{\alpha-2}$; such distributions are remainders in the present analysis. We use first the paracontrolled structure of $u$ in the term $\Delta u$, and the continuity result on the commutator $\big[\Delta,\overline\Pi\,\big]$, given in Lemma 5.1 of \cite{GIP}, to write
\begin{equation*}
\begin{split}
\big(a(u)-a(u_0^T)\big)\Delta u &= -\big(a(u)-a(u_0^T)\big)\overline{\Pi}_{u'}\xi + (\Asterisk) + \big(a(u)-a(u_0^T)\big)\Delta u^\sharp   \\
													 &= -\Pi_{(a(u)-a(u_0^T))u'} \xi - \Pi\Big(a(u)-a(u_0^T), \overline\Pi_{u'}\xi\Big) + (\Asterisk)   \\
													 &\quad\,+ \big(a(u)-a(u_0^T)\big)\Delta u^\sharp;
\end{split}
\end{equation*}
the second equality is a special case of Theorem 6 in \cite{BailleulBernicot2}. Note that the term in $\Delta u^\sharp$ cannot go inside the remainder as it has an explosive spatial $C^{2\beta}$ norm at time $0^+$. Using that $\big(\overline\Pi_{u'}\xi - \Pi_{u'}\xi\big)\in C_TC^{\alpha+\beta-2}$, such as proved in Lemma 5.1 of \cite{GIP}, the first continuity estimate of Lemma \ref{LemFixedPointSetting} then gives
\begin{equation}
\label{EqRenormalisation2}
\Pi\Big(a(u), \overline\Pi_{u'}\xi\Big) = a'(u)\,(u')^2\,\Pi(X,\xi) + a'(u)u'\Pi\big(u^\sharp,\xi\big) + (\Asterisk);
\end{equation}
we have in addition
$$
\Big\|\Pi\big(a(u_0^T),\overline\Pi_{u'}\xi\big)\Big\|_{C^{\alpha+2\beta-2}} \lesssim \|a\|_{\mcC^1}\left(1+ T^{-\frac{2\beta-\alpha}{2}}\|u_0\|_{C^\alpha}\right)\|u'\|_{\mcC^\beta} \|\xi\|_{C^{\alpha-2}}. 
$$
(Schauder estimates for the resolution operator of $\mathscr{L}^0$ will later take care of the exploding factor in $T$.) We can thus rewrite equation \eqref{EqRewritegPAM} at this point under the form
$$
\mathscr{L}^0 u = \Pi_{g(u)-(a(u)-a(u_0^T))u'}\xi + \big(a(u)-a(u_0^T)\big)\Delta u^\sharp + (\Asterisk_1);
$$
building on the second identity of Lemma \ref{LemFixedPointSetting}, we end up with the equation
$$
\Pi_{a(u_0^T)u'}\xi + \mathscr{L}^0 u^\sharp = \Pi_{g(u)-(a(u)-a(u_0^T))u'}\xi + \big(a(u)-a(u_0^T)\big)\Delta u^\sharp + (\Asterisk_2),
$$
that is
$$
\mathscr{L}^0 u^\sharp = \Pi_{g(u)-a(u)u'}\xi + \big(a(u)-a(u_0^T)\big)\Delta u^\sharp + (\Asterisk_2).
$$
We use here the notation $(\Asterisk_i)$ to single out these particular terms, as opposed to the above unspecified remainders; each of them takes the form 
$$
(\Asterisk_i) =  \big(g'(u)-a'(u)u'\big)\,\Pi\big(u^\sharp,\xi\big) + (\Asterisk_i)_{\geq \alpha+\beta-2},
$$
with $(\Asterisk_i)_{\geq \alpha+\beta-2}\in C_TC^{\geq \alpha+\beta-2}$. As said above, the term involving $u^\sharp$ needs to be evaluated in a weighted space, although it has positive space regularity at each positive time. Let say that a \textit{constant depends on the data} if it depends on $\|\xi\|_{C^{\alpha-2}}$, $\|X\|_{C^\alpha}$, $\big\|\Pi(X,\xi)\big\|_{C^{2\alpha-2}}$, $\|g\|_{C^3_b}$, $\|a\|_{C^3_b}$ and possibly $\|u_0\|_{C^\alpha}$. Given the fact that $(\Asterisk_i)_{\geq \alpha+\beta-2}$ is given explicitly in terms of multilinear functions of $u,u'$ or $C^2_b$ functions of $u$, the term $(\Asterisk_i)_{\geq \alpha+\beta-2}$ defines a function of $(u,u')\in {\sf C}_{\alpha,\beta}(X)$ that is locally Lipschitz, with a Lipschitz constant that depends polynomially on $\big\|(u,u')\big\|_{\alpha,\beta}$, and $(\Asterisk_i)_{\geq \alpha+\beta-2}$ has itself a $C_TC^{\geq \alpha+\beta-2}$-norm that is polynomial in terms of $\big\|(u,u')\big\|_{\alpha,\beta}$; everything depends of course on the data.

\medskip

For $(u,u')\in{\sf C}_{\alpha,\beta}(X)$, set 
$$
\Phi(u,u') := (v,v'),
$$
where
\begin{equation*}
\begin{split}
v' &:= \frac{g(u)-\big(a(u)-a(u_0^T)\big)u'}{a(u_0^T)}   \\
\mathscr{L}^0 v &:= \Pi_{a(u_0^T)v'}\xi + \big(a(u)-a(u_0^T)\big)\Delta u^\sharp + (\Asterisk_1),
\end{split}
\end{equation*}
with $v_{t=0} = u_0$, and 
\begin{equation}
\label{EqDynamicsV}
\mathscr{L}^0 v^\sharp := \big(a(u)-a(u_0^T)\big)\Delta u^\sharp + (\Asterisk_2),
\end{equation}
with $v^\sharp_{t=0} = u_0-\overline\Pi_{v'_{t=0}}X$. For $\lambda$ positive, set 
$$
\frak{B}_T(\lambda) := \left\{(u,u')\in {\sf C}_{\alpha,\beta}(X)\,;\, u_{t=0}=u_0,\;u'_{t=0} = \frac{g(u_0)}{a(u_0)},\; \big\|(u,u')\big\|_{\alpha,\beta} \leq \lambda \right\}.
$$
We are going to prove that \vspace{0.1cm}

\begin{itemize}
   \item[\textcolor{gray}{$\bullet$}] the map $\Phi$ sends $\frak{B}_T(\lambda)$ into itself, for an adequate choice of radius $\lambda$ and a choice of sufficiently small time horizon $T$,   \vspace{0.1cm}
   
   \item[\textcolor{gray}{$\bullet$}] it is in that case a contraction of $\frak{B}_T(\lambda)$.
\end{itemize}

\bigskip

\subsection{Fixed point}
\label{SubsectionFixedPoint}

We first give in the next lemma a control on the  two terms involving $v^\sharp$ in the dynamics \eqref{EqDynamicsV} of $v$.

\medskip

\begin{lem}
\label{LemVSharp}
Given ${\sf u} := (u,u'),\,{\sf u}_1:=(u_1,u'_2),\,{\sf u}_2:=(u_2,u_2')\in \frak{B}_T(\lambda)$, set
\begin{equation*}
\begin{split}
&\varep_1({\sf u}) := \varep_1(u,u') = \big(g'(u)-a'(u)u'\big) \, \Pi(u^\sharp,\xi),   \\
&\varep_2({\sf u}) := \varep_2(u,u') = \big(a(u)-a(u_0^T)\big) \Delta u^\sharp.
\end{split}
\end{equation*}
Then we have the estimates
\begin{equation*}
\begin{split}
&\underset{0< t\leq T}{\sup}\; t^{\frac{2\beta-\alpha}{2}}\,\big\|\varep_1({\sf u})(t)\big\|_{C^{2\alpha-2}} \lesssim C_1\big(\|{\sf u}\|_{\alpha,\beta}\big),   \\
&\underset{0< t\leq T}{\sup}\; t^{\frac{2\beta-\alpha}{2}}\,\big\|\varep_2({\sf u})(t)\big\|_{C^{2\beta-2}} \lesssim T^\frac{\alpha-\beta}{2}\, C_1\big(\|{\sf u}\|_{\alpha,\beta}\big) + C_2
\end{split}
\end{equation*}
for a constant $C_1\big(\|{\sf u}\|_{\alpha,\beta}\big)$ depending polynomially on the data and $\|{\sf u}\|_{\alpha,\beta}$, and $C_2$ depending on the data, and 
{\small 
\begin{equation*}
\begin{split}
&\underset{0< t\leq T}{\sup}\; t^{\frac{2\beta-\alpha}{2}}\,\Big\|\big(\varep_i({\sf u}_1) - \varep_i({\sf u}_2)\big)(t)\Big\|_{C^{2\alpha-2}} \lesssim C_3\big(\|{\sf u}_1\|_{\alpha,\beta}, \|{\sf u}_2\|_{\alpha,\beta}\big)\,\big\|{\sf u}_1 - {\sf u}_2\big\|_{\alpha,\beta},   \\
&\underset{0< t\leq T}{\sup}\; t^{\frac{2\beta-\alpha}{2}}\,\Big\|\big(\varep_i({\sf u}_1) - \varep_i({\sf u}_2)\big)(t)\Big\|_{C^{2\beta-2}} \lesssim T^{\frac{\alpha-\beta}{2}}\,C_3\big(\|{\sf u}_1\|_{\alpha,\beta}, \|{\sf u}_2\|_{\alpha,\beta}\big)\,\big\|{\sf u}_1 - {\sf u}_2\big\|_{\alpha,\beta},
\end{split}
\end{equation*}   }
for a constant $C_3\big(\|{\sf u}_1\|_{\alpha,\beta},\|{\sf u}_2\|_{\alpha,\beta}\big)$ depending polynomially on the data and $\|{\sf u}_1\|_{\alpha,\beta} $ and $\|{\sf u}_2\|_{\alpha,\beta} $.
\end{lem}

\medskip

Remark the gain of a factor $T^{\frac{\alpha-\beta}{2}}$ in the estimate for the local Lipschitz character of $\varep_2$ as a function of ${\sf u}$; this is taken care of by Schauder estimates for $\varepsilon_1$.

\medskip

\begin{Dem}
The size bound for $\varepsilon_i({\sf u})(t)$ is elementary. For $\varepsilon_1({\sf u})$, write 
\begin{equation*}
\begin{split}
\Big\| \big\{\big(g' (u)-a'(u)u'\big)\Pi\big(u^\sharp, \xi\big)\big\}(t)\Big\|_{C^{2 \alpha - 2}} &\lesssim \big\|\big(g' (u)-a'(u)u'\big)(t)\big\|_{C^\beta}\,\big\| u^\sharp(t)\big\|_{C^{2 \beta}} \, \| \xi \|_{\alpha - 2}   \\
&\lesssim (\cdots)\, \big\| u^\sharp(t)\big\|_{C^{2 \beta}} \, \| \xi \|_{C^{\alpha - 2}}
\end{split}
\end{equation*}
with 
$$
(\cdots) = \| g'\|_{C^1} \big((1 + \| u \|_{\mcC^\alpha}\big) + \| a'\|_{C^1} \big(1 + \| u \|_{\mcC^\alpha}\big)\,\| u' \|_{\mcC^\beta},
$$
and use the fact that $3\beta-2$ is positive to get
\begin{equation*}
\begin{split}
\big\| \big\{\big(a(u)-a(u_0^T)\big)\Delta u^\sharp\big\}(t) \big\|_{C^{2 \beta - 2}} &\lesssim \big\| a (u) - a(u^T_0) \big\|_{C^\beta} \, \big\| u^\sharp(t)\big\|_{C^{2 \beta}} .
\end{split}
\end{equation*}
This can be further estimated using
$$
\big\| a(u)-a(u_0^T) \big\|_{\mcC^\beta} \leq \big\| a (u) - a (u_0) \big\|_{\mcC^\beta} + \big\| a (u_0) - a (u_0^T) \big\|_{C^\beta},
$$
where
$$
\big\| a (u) - a (u_0) \big\|_{\mcC^\beta} \lesssim T^{\frac{\alpha - \beta}{2}} \| a \|_{C^2} \Big(1 + \| u \|_{\mcC^\alpha}\Big) \| u - u_0\|_{\mcC^\alpha},
$$
and  from Lemma \ref{lem:aux} in the Appendix
{\small \begin{equation*}
\begin{split}
\big\| a (u_0) - a (u_0^T) \big\|_{C^\beta} &\leq \| a' \|_{C^1} \big\| u_0 - u_0^T \big\|_{L^{\infty}} + \| a'' \|_{C^0} \big\| u_0 - u_0^T \big\|_{L^{\infty}} \| u_0\|_{C^\beta} + \| a' \|_{C^0} \| u_0 - u_0^T \|_{C^\beta}   \\
&\lesssim T^{\frac{\alpha - \beta}{2}} \| u_0 \|_{C^\alpha}.
\end{split}
\end{equation*}   }
Therefore
\begin{equation*}
\begin{split}
\big\| \big\{\big(a(u)-a(u_0^T)\big)\Delta u^\sharp\big\}(t) \big\|_{C^{2 \beta - 2}} &\lesssim T^{\frac{\alpha - \beta}{2}}  \big(1 + \| u \|_{\mcC^\alpha}+\|u_0\|_{C^\alpha}\big) \, \big\| u^\sharp (t)\big\|_{C^{2 \beta}}.
\end{split}
\end{equation*}
We only look at the Lipschitz estimate for $\varep_2$ and leave the  reader treat the easier case of $\varepsilon_1$. It suffices in the former case to write
{\small \begin{equation*}
\begin{split}
&\Big\| \Big\{\big(a (u_1) - a (u^T_0)\big) \Delta u_1^\sharp - \big(a (u_2) - a (u^T_0)\big) \Delta u_2^\sharp\Big\}(t)\Big\|_{C^{2 \beta - 2}}   \\
&\leq \Big\| \big\{\big(a (u_1) - a (u_2)\big) \Delta u_1^\sharp\big\}(t) \Big\|_{C^{2 \beta - 2}} + \Big\| \Big\{\big(a (u_2) - a (u^T_0)\big) \Delta (u_1^\sharp - u_2^\sharp)\Big\}(t) \Big\|_{C^{2 \beta - 2}}   \\
&\lesssim \big\| \big(a (u_1) - a (u_2)\big(t) \big\|_{C^\beta} \big\| u_1^\sharp\big\|_{C^{2 \beta}} + \Big(\big\| a (u_1) - a (u_0) \big\|_{C^\beta} + \big\| a (u_0) - a (u_0^T) \big\|_{C^\beta}\Big) \, \big\| u_1^\sharp - u_2^\sharp \big\|_{C^{2 \beta}}   \\
&\lesssim T^{\frac{\alpha - \beta}{2}} \| a \|_{C^2} \big(1 + \| u_1\|_{\mcC^\alpha}\big) \, \| u_1 - u_2 \|_{\mcC^\alpha} \big\| u_1^\sharp(t)\big\|_{C^{2 \beta}}   \\
&\quad+ T^{\frac{\alpha - \beta}{2}} \Big(\| a \|_{C^2} \big(1 + \| u_1\|_{\mcC^\alpha}\big) + \| u_0 \|_{C^\alpha}\Big) \big\| u_1^\sharp - u_2^\sharp \big\|_{C^{2\beta}}
\end{split}
\end{equation*}   }
to get the result.
\end{Dem}

\medskip

\begin{lem}
\label{LemmaV'}
For ${\sf u}_1, {\sf u}_2$ in $\frak{B}_T(\lambda)$ and ${\sf v}_1 := \Phi({\sf u}_1) = \big(v_1,v_1'\big)$ and ${\sf v}_2 := \Phi({\sf u}_2) = \big(v_2,v'_2\big)$, we have the estimates
\begin{equation*}
\begin{split}
&\big\| v_1' \|_{\mcC^\beta} \lesssim T^{\frac{\alpha - \beta}{2}} C_1\big(\|{\sf u}_1\|_{\alpha,\beta}\big) + C_2,   \\
&\| v_2' - v_1' \|_{\mcC^\beta} \lesssim T^{\frac{\alpha - \beta}{2}} \, C_3\big(\|{\sf u}_1\|_{\alpha,\beta}, \|{\sf u}_2\|_{\alpha,\beta}\big) \, \big\| {\sf u}_1 - {\sf u}_2 \big\|_{\alpha,\beta}.
\end{split}
\end{equation*}
\end{lem}

\medskip

\begin{Dem}
We first bound $\big\| v_1' \|_{\mcC^\beta}$ and start for that purpose from the rough estimate
$$
\| v_1' \|_{\mcC^\beta} \leq \left\| \frac{1}{a(u_0^T)} \right\|_{C^\beta} \Big(\big\| g (u_1) \big\|_{\mcC^\beta} + \big\| \big(a(u_1)-a(u_0^T)\big) u_1' \big\|_{\mcC^\beta}\Big).
$$
Lemma \ref{lem:aux} gives on the one hand
$$
\big\| g (u_1) \big\|_{\mcC^\beta} \lesssim T^{\frac{\alpha - \beta}{2}} \| g \|_{C^1} \Big(1 + \| u_1 \|_{\mcC^\alpha}\Big) + \big\| g (u_0)\big\|_{C^\beta}.
$$
(Be careful that $u_1$ is \textit{not} the time $1$ value of some $u$.) On the other hand, as in the proof of Lemma \ref{LemVSharp}, we have
\begin{equation*}
\label{eq:b}
\Big\| \big(a(u_1)-a(u_0^T)\big) u_1'\Big\|_{\mcC^\beta} \lesssim T^\frac{\alpha-\beta}{2}\big(1+\|u_1\|_{C^\alpha}+\|u_0\|_{C^\alpha}\big)\| u_1' \|_{\mcC^\beta} .
\end{equation*}

\medskip
  
In order to obtain the Lipschitz bound of the statement, Lemma \ref{lem:aux} gives us again
$$
\big\| g (u_1) - g (u_2) \big\|_{\mcC^\beta} \lesssim T^{\frac{\alpha -\beta}{2}} \| g \|_{C^2} \Big(1 + \| u_1 \|_{\mcC^\alpha}\Big) \| u_1 - u_2\|_{\mcC^\alpha},
$$
and, with $b_i := a(u_i)-a(u_0^T)$, for $i=1,2$,
\begin{equation*}
\begin{split}
\big\| b_1 u_1' - b_2u_2' \big\|_{\mcC^\beta} &\leq \|b_1\|_{\mcC^\beta} \| u_1' - u_2' \|_{C_T L^{\infty}}+\|b_1\|_{C_T L^{\infty}} \| u_1' - u_2' \|_{\mcC^\beta} \\
&\quad+ \|b_1-b_2\|_{\mcC^\beta}\|u'_2\|_{C_TL^\infty}+ \| b_1-b_2 \|_{C_T L^{\infty}} \| u_2'\|_{\mcC^\beta}   \\
&\lesssim T^{\frac{\beta}{2}} \| a \|_{C^1} \Big(1 + \| u_1\|_{\mcC^\alpha} + \| u_0 \|_{\alpha}\Big) \| u_1' - u_2'\|_{\mcC^\beta}   \\
&\quad+ T^\frac{\alpha}{2}\|a\|_{C^1}\|u_1\|_{\mathcal{C}^\alpha}\| u_1' - u_2' \|_{\mcC^\beta}\\
&\quad+ T^\frac{\alpha-\beta}{2}\|a\|_{C^2}\big(1+\|u_1\|_{\mcC^\alpha}\big)\|u_1-u_2\|_{\mcC^\alpha}\|u_2'\|_{\mcC^\alpha}\\
&\quad+ T^{\frac{\alpha}{2}} \| a \|_{C^2} \big(1 + \| u_2\|_{\mcC^\alpha}\big) \| u_1 - u_2 \|_{\mcC^\alpha} \| u_2'\|_{\mcC^\beta}.
\end{split}
\end{equation*} 
\end{Dem}

\medskip

\begin{lem}
\label{lem:a+b}
Let an initial condition $f_0 \in C^\alpha$ be given, together with another function $g \in C^2$, bounded below by a positive constant. Let also $\phi_1\in C\big((0,T],C^{2\beta-2}\big)$ with 
\begin{equation}
\label{EqConditionPhi1}
\underset{0<t\leq T}{\sup}\, t^\frac{2\beta-\alpha}{2} \, \big\|\phi_1(t)\big\|_{C^{2 \beta-2}} < \infty,
\end{equation}
and 
$\phi_2\in C\big((0,T],C^{\alpha+\beta-2}\big)$ with 
\begin{equation}
\label{EqConditionPhi2}
\underset{0<t\leq T}{\sup}\, t^\frac{2\beta-\alpha}{2} \, \big\|\phi_2(t)\big\|_{C^{\alpha+\beta-2}} < \infty
\end{equation}
be given. Let $f$ be the solution to the evolution equation
\begin{equation}
\label{Eqf}
\partial_t f - g \Delta f = \phi_1 + \phi_2, \qquad f (0) = f_0.
\end{equation}
Then, choosing the time horizon $T$ small enough, we have the estimate
\begin{equation}
\label{EqEstimate}
\begin{split}
\underset{0<t\leq T}{\sup} \, &t^\frac{2\beta-\alpha}{2} \, \big\| f (t) \big\|_{C^{2 \beta}} + \| f\|_{\mcC^\alpha}   \\
&\lesssim \| f_0 \|_{C^\alpha} + \underset{0<t\leq T}{\sup} \, t^\frac{2\beta-\alpha}{2} \, \big\| \phi_1(t) \big\|_{C^{2 \beta - 2}} + T^{\frac{\alpha - \beta}{2}} \, \underset{0<t\leq T}{\sup} \, t^\frac{2\beta-\alpha}{2} \, \big\| \phi_2(t) \big\|_{C^{\alpha +\beta - 2}},
\end{split}
\end{equation}
with an implicit multiplicative positive constant in the right hand side depending only on the $C^\alpha$-norm of $g$.
\end{lem}

\medskip

The fact that this multiplicative positive constant depends only on the $C^\alpha$-norm of $g$ is crucial for what comes next.

\medskip

\begin{Dem}
Let $(Q_t)_{0\leq t\leq T}$ stand for the semigroup generated by the operator $\textrm{div}(g \nabla\cdummy)$. We know from \cite{AuscherMcIntosh} and \cite{BailleulBernicot1} that the resolution operator associated with the heat operator built from $\textrm{div}(g \nabla\cdummy)$ satisfies the Schauder estimates, such as stated in Lemma A.7-A.9 of \cite{GIP} and Corollary 4.5 in \cite{cannizzaro_multidimensional_2015}, with implicit multiplicative constants depending only on the $C^\alpha$-norm of $g$. Write $\nabla g \cdot \nabla f$ for $\sum_i \partial_i g\,\partial_i f$. The solution $f$ to equation \eqref{Eqf} is given in mild formula   
$$
f_t = Q_t f_0 - \int_0^t Q_{t - s} \big(\nabla g \cdummy \nabla f(s)\big) \,\mathd s + \int_0^t Q_{t - s} \phi_1 (s) \,\mathd s + \int_0^t Q_{t - s} \phi_2 (s) \, \mathd s.
$$
Note that since the exponent $(2 \beta + \alpha - 2)$ is positive, we have
$$
\underset{0<t\leq T}{\sup} \, t^\frac{2\beta-\alpha}{2} \, \big\| \nabla g \cdummy \nabla f(t) \big\|_{C^{\alpha - 1}} \lesssim \| g \|_{C^\alpha} \underset{0<t\leq T}{\sup} \, t^\frac{2\beta-\alpha}{2} \, \big\| f (t)\big\|_{C^{2 \beta}}.
$$
It follows from the Schauder estimates that we have at any positive time $t$ in $(0,T]$ the upper bound
\begin{equation*}
\begin{split}
t^\frac{2\beta-\alpha}{2} \, \big\| f(t) \big\|_{C^{2 \beta}} &\lesssim \underset{0<s\leq t}{\sup} \, s^\frac{2\beta-\alpha}{2} \, \big\| \phi_1(s)\big\|_{C^{2 \beta - 2}} + T^{\frac{\alpha - \beta}{2}} \, \underset{0<s\leq t}{\sup} \, s^\frac{2\beta-\alpha}{2} \, \big\| \phi_2(s) \big\|_{C^{\alpha +\beta - 2}}   \\
&\quad+ T^{\frac{1+\alpha-2\beta}{2}} \,\| g \|_{C^\alpha} \underset{0<s\leq t}{\sup} \, s^\frac{2\beta-\alpha}{2} \,\big\| f(s) \big\|_{C^{2 \beta}} + \| f_0 \|_{C^\alpha};
\end{split}
\end{equation*}
taking the time horizon $T$ small enough then yields part of the estimate of the statement. Next, we have
\begin{equation*}
\begin{split}
\big\| f(t) - f(s) \big\|_{L^{\infty}} \leq &\big\| (Q_{t - s} - \tmop{Id}) f_0\big\|_{L^{\infty}} + \left\| \int_s^t Q_{t - r} (\nabla g \cdummy \nabla f) \,\mathd r\right\|_{L^{\infty}}   \\
&+ \left\| \int_0^s (Q_{t - s} - \tmop{Id}) Q_{s -r} (\nabla g \cdummy \nabla f) \mathd r \right\|_{L^{\infty}} + \left\| \int_s^t Q_{t - r} \phi_1(r)\, \mathd r \right\|_{L^{\infty}}   \\
&+ \left\| \int_0^s (Q_{t - s} - \tmop{Id}) Q_{s - r} \phi_1(r) \mathd r\right\|_{L^{\infty}} + \left\| \int_s^t Q_{t - r} \phi_2(r) \,\mathd r \right\|_{L^{\infty}}   \\
&+ \left\| \int_0^s (Q_{t - s} - \tmop{Id}) Q_{s - r} \phi_2(r)\, \mathd r\right\|_{L^{\infty}} = : {\sf I}_1 + \cdots + {\sf I}_7.
\end{split}
\end{equation*}
Using Lemma A.8 {\cite{GIP}} to the first term, we obtain
$$
{\sf I}_1 \lesssim | t - s |^\frac{\alpha}{2} \, \| f_0 \|_{C^\alpha}.
$$
For the other terms, and for any positive exponent $a$, we use repeatedly the following elementary extension of Lemma A.7 of \cite{GIP}
\begin{equation}\label{eq:sch}
\big\|Q_tu\big\|_{L^\infty} \lesssim t^{-\frac{a}{2}} \, \|u\|_{C^{-a}}.
\end{equation}
(It can be seen to hold as follows. Writing $L$ for the operator $\textrm{div}(g \nabla\cdummy)$ and setting $R_s := (sL)e^{-sL}$, we know that $R_su$ is bounded in $L^\infty$ by $s^{-a/2}$ if $u$ is $C^{-a}$ -- this semigroup picture of H\"older spaces is explained and used for instance in \cite{BailleulBernicot1}. The above continuity estimate comes then from the integral representation $Q_t = \int_t^\infty R_s\,\frac{ds}{s}$.) Apply then \eqref{eq:sch} to the second term and Lemma A.8 of \cite{GIP} to the third one to get
\begin{equation*}
\begin{split}
{\sf I}_2 + {\sf I}_3 &\lesssim \int_s^t (t - r)^{\frac{\alpha - 1}{2}} \big\| \nabla g \cdummy \nabla f(r) \big\|_{C^{\alpha - 1}} \, \mathd r + | t - s |^\frac{\alpha}{2} \, \int_0^s \Big\| Q_{s - r} \big(\nabla g \cdummy \nabla f(r)\big) \Big\|_{C^\alpha} \, \mathd r   \\
&\lesssim \| g \|_{C^\alpha} \left\{\underset{0<r\leq t}{\sup}\, r^\frac{2\beta-\alpha}{2} \, \big\| f(r) \big\|_{C^{2 \beta}}\right\} \int_s^t (t - r)^{\frac{\alpha - 1}{2}}\, r^{-\frac{2 \beta - \alpha}{2}} \, \mathd r   \\
&\quad+ | t - s |^\frac{\alpha}{2} \| g \|_{C^\alpha} \left\{\underset{0<r\leq t}{\sup}\, r^\frac{2\beta-\alpha}{2} \, \big\| f(r) \big\|_{C^{2 \beta}}\right\} \int_0^s ( s - r)^{- \frac{1}{2}} \, r^{-\frac{2\beta-\alpha}{2}} \, \mathd r   \\
&\lesssim \Big(| t - s |^{\frac{1}{2} + \alpha - \beta} + T^{\frac{1}{2} - \beta + \frac{\alpha}{2}} \, | t - s |^\frac{\alpha}{2}\Big) \, \| g \|_{C^\alpha} \left\{\underset{0<r\leq t}{\sup}\, r^\frac{2\beta-\alpha}{2} \, \big\| f(r) \big\|_{C^{2 \beta}}\right\}   \\
&\lesssim T^\frac{1+\alpha-2\beta}{2}\, | t - s |^\frac{\alpha}{2}\, \| g \|_{C^\alpha} \left\{\underset{0<r\leq t}{\sup}\, r^\frac{2\beta-\alpha}{2} \, \big\| f(r) \big\|_{C^{2 \beta}}\right\}.
\end{split}
\end{equation*}
We bound similarly the quantities
\begin{equation*}
\begin{split}
{\sf I}_4 + {\sf I}_5 &\lesssim \int_s^t (t - r)^{- 1 + \beta} \, \big\| \phi_1(r) \big\|_{C^{2\beta - 2}} \mathd r + | t - s |^\frac{\alpha}{2} \int_0^s \big\| Q_{s - r} \phi_1(r) \big\|_{C^\alpha} \,\mathd r   \\
&\lesssim | t - s |^\frac{\alpha}{2} \left\{\underset{0<r\leq t}{\sup}\, r^\frac{2\beta-\alpha}{2}\,\big\| \phi_1(r) \big\|_{C^{2 \beta - 2}}\right\} \int_s^t (t - r)^{- 1 +\beta - \frac{\alpha}{2}}\, r^{-\frac{2\beta-\alpha}{2}} \,\mathd r   \\
&\quad+ | t - s |^\frac{\alpha}{2} \left\{\underset{0<r\leq t}{\sup}\, r^\frac{2\beta-\alpha}{2}\, \big\| \phi_1(r) \big\|_{C^{2\beta - 2}}\right\} \int_0^s (s - r)^{- 1 + \beta  - \frac{\alpha}{2}}\, r^{-\frac{2\beta-\alpha}{2}}\, \mathd r
\end{split}
\end{equation*}
and, since $2\beta-\alpha \leq \beta$, and one can assume $0\leq t\leq T\leq 1$, 
\begin{equation*}
\begin{split}
{\sf I}_6 + {\sf I}_7 &\lesssim \int_s^t (t - r)^{- 1 + \frac{\alpha}{2} + \frac{\beta}{2}} \big\| \phi_2(r) \big\|_{C^{\alpha + \beta - 2}}\, \mathd r + | t - s |^\frac{\alpha}{2} \int_0^s \big\| Q_{s - r} \phi_2(r)\big\|_{C^\alpha} \,\mathd r   \\
&\lesssim | t - s |^\frac{\alpha}{2}\, T^\frac{\alpha-\beta}{2} \, \left\{\underset{0<r\leq t}{\sup}\, r^\frac{2\beta-\alpha}{2} \big\|\phi_2(r)\big\|_{C^{\alpha + \beta - 2}}\right\} \int_s^t (t - r)^{- 1 + \frac{\beta}{2}}\, r^{- \frac{\beta}{2}}\, \mathd r   \\
&\quad+ | t - s |^\frac{\alpha}{2}\, T^\frac{\alpha-\beta}{2} \,  \left\{\underset{0<r\leq t}{\sup}\, r^\frac{2\beta-\alpha}{2} \big\| \phi_2(r) \big\|_{C^{\alpha + \beta - 2}}\right\} \int_0^s (s - r)^{- 1 + \frac{\beta}{2}} \,r^{- \frac{\beta}{2}}\, \mathd r.
\end{split}
\end{equation*}
Since for any fixed positive exponent $\delta\in (0,1)$, we have
$$
\int_0^t (t - r)^{- 1 + \delta}\, r^{- \delta} \,\mathd r \lesssim 1,
$$
uniformly in $t\in (0,T]$ and $T\leq 1$, we deduce that
\begin{equation*}
\begin{split}
\| f\|_{C^{\alpha / 2}_T L^{\infty}} &\lesssim \| f_0 \|_{C^\alpha} + \underset{0<t\leq T}{\sup}\, t^\frac{2\beta-\alpha}{2}\, \big\| \phi_1(t) \big\|_{C^{2 \beta - 2}}   \\
&\quad+ T^{\frac{\alpha - \beta}{2}} \underset{0<t\leq T}{\sup}\, t^\frac{2\beta-\alpha}{2} \,\big\| \phi_2(t) \big\|_{C^{\alpha + \beta - 2}}.
\end{split}
\end{equation*}
Very similar arguments give the estimate
\begin{equation*}
\begin{split}
\| f_t \|_{C^\alpha} &\lesssim \| f_0 \|_{\alpha} + \| g \|_{C^\alpha} \underset{0<r\leq t}{\sup}\, r^\frac{2\beta-\alpha}{2} \big\| f(r)\big\|_{C^{2 \beta}} \int_0^t (t - r)^{-\frac{1}{2}} \,r^{- \frac{2\beta-\alpha}{2} }\,\mathd r   \\
&\quad+ \underset{0<r\leq t}{\sup}\, r^\frac{2\beta-\alpha}{2} \big\| \phi_1(r) \big\|_{C^{2 \beta - 2}} \int_0^t (t - r)^{- 1 + \beta - \frac{\alpha}{2}} \, r^{- \frac{2\beta-\alpha}{2}} \,\mathd r   \\
&\quad+ \underset{0<r\leq t}{\sup}\, r^\frac{\beta}{2} \big\| \phi_2(r) \big\|_{C^{\alpha + \beta - 2}} \int_0^t (t - r)^{- 1 + \frac{
\beta}{2}} \, r^{- \frac{\beta}{2}}\, \mathd r,
\end{split}
\end{equation*}     
from which we finally get
\begin{equation*}
\begin{split}
\| f \|_{L^\infty_T \mathcal{C}^{\alpha}} &\lesssim \| f_0 \|_{C^\alpha} + \underset{0<t\leq T}{\sup}\,  t^\frac{2\beta-\alpha}{2} \,\big\| \phi_1(t) \big\|_{C^{2 \beta - 2}}   \\
 &\quad+ T^{\frac{\alpha - \beta}{2}} \underset{0<t\leq T}{\sup}\, t^\frac{2\beta-\alpha}{2} \,\big\| \phi_2(t) \big\|_{C^{\alpha + \beta - 2}},
\end{split}
\end{equation*} 
and the result of the statement.
\end{Dem}

\medskip

Recall that $u_0$ is $\alpha$-H\"older. The following statement is a direct corollary of Lemma~\ref{lem:a+b} and the fact that while $u_0^T$ is regular as a consequence of the regularizing properties of the heat semigroup, its norm as a regular element blows up as $T$ decreases to $0$. On the other hand, the spatial $\alpha$-H\"older norm of $u_0^T$, or $a(u_0^T)$, is controlled in terms of the $\alpha$-H\"older norm of $u_0$, with no exploding factor, uniformly in $T$ near $0^+$.

\medskip

\begin{cor}
\label{Corollary}
Assume we are given some functions $\phi_1\in C\big((0,T],C^{2\beta-2}\big)$ and $\phi_2\in C\big((0,T],C^{\alpha+\beta-2}\big)$ satisfying the estimates \eqref{EqConditionPhi1} and \eqref{EqConditionPhi2}. Given $z_0\in C^\alpha$ and $z'\in\mcC^\beta$, let $z$ stand for the solution of the quation 
\begin{equation*}
\label{eq:parac}
\big(\partial_t - a(u_0^T)\Delta\big) z = \Pi_{a (u_0^T) z'} \xi + \phi_1 + \phi_2, \qquad z (0) = z_0.
\end{equation*}
Then $(z, z') \in {\sf C}_{\alpha, \beta}(X)$, and we have the size estimate
\begin{equation}
\label{EqEstimateParac}
\begin{split}
\underset{0<t\leq T}{\sup}\,t^\frac{2\beta-\alpha}{2}\, &\big\| z^\sharp(t)\big\|_{C^{2 \beta}} + \big\| z^\sharp\big\|_{\mcC^\alpha}   \\
&\lesssim \| z_0 \|_{C^\alpha} + \big\| z'(0)\big\|_{L^{\infty}} \| X \|_{C^\alpha} + T^{\frac{\alpha - \beta}{2}} \big(1 + \| u_0\|_{C^\alpha}\big) \| z' \|_{\mcC^\beta} \| X \|_{C^\alpha}   \\
&\quad \underset{0<t\leq T}{\sup}\,t^\frac{2\beta-\alpha}{2}\,\big\|\phi_1(t)\big\|_{C^{2 \beta - 2}} + T^{\frac{\alpha - \beta}{2}} \underset{0<t\leq T}{\sup}\,t^\frac{2\beta-\alpha}{2} \big\|\phi_2(t)\big\|_{C^{\alpha + \beta - 2}},
\end{split}
\end{equation}   
with an implicit multiplicative positive constant in the right hand side depending only on the $C^\alpha$-norm of $u_0$. If $(y,y')\in {\sf C}_{\alpha,\beta}(X)$ is associated similarly to another set of data $\psi_1,\psi_2, y_0$ and $y'$, with $y$ solution of the equation
$$
\big(\partial_t - a(u_0^T)\Delta\big) y = \Pi_{a (u_0^T) y'} (\xi) + \psi_1 + \psi_2, \qquad y (0) = y_0 \in C^\alpha,
$$
then
\begin{equation*}
\begin{split}
&\underset{0<t\leq T}{\sup}\,t^\frac{2\beta-\alpha}{2}\, \big\| z^\sharp(t)-y^\sharp(t)\big\|_{C^{2 \beta}} + \big\| z^\sharp-y^\sharp\big\|_{\mcC^\alpha}   \\
&\lesssim \| z_0 - y_0\|_{C^\alpha} + \big\| z' (0) - y' (0) \big\|_{L^{\infty}} \| X \|_{C^\alpha} + T^{\frac{\alpha - \beta}{2}} \big(1 + \| u_0 \|_{C^\alpha}\big) \| z' - y' \|_{\mcC^\beta} \| X\|_{C^\alpha}   \\
&\quad+ \underset{0<t\leq T}{\sup}\,t^\frac{2\beta-\alpha}{2}\,\big\| \phi_1(t) - \psi_1(t)\big\|_{C^{2 \beta - 2}} + T^{\frac{\alpha - \beta}{2}} \underset{0<t\leq T}{\sup}\,t^\frac{2\beta-\alpha}{2}\,\big\|\phi_2(t) - \psi_2(t)\big\|_{C^{\alpha + \beta - 2}},
\end{split}
\end{equation*}     
here again, with an implicit multiplicative positive constant in the right hand side depending only on the $C^\alpha$-norm of $u_0$.
\end{cor}

\medskip

\begin{Dem}
Recall we write $\mathscr{L}^0$ for the operator $\partial_t - a(u_0^T)\Delta$, for short. Set 
$$
z^\sharp := z - \overline{\Pi}_{z'}X.
$$
As this function is the solution of the equation
\begin{equation*}
\begin{split}
\mathscr{L}^0z^\sharp &= \mathscr{L}^0 z - \mathscr{L}^0\Big(\overline\Pi_{z'}X\Big)   \\
&= \Pi_{a (u_0^T) z'} \xi + \phi_1 + \phi_2 - \Big\{\mathscr{L}^0 \Big(\overline{\Pi}_{z'} X\Big) + \Pi_{a (u_0^T) z'} (-\xi)\Big\} + \Pi_{a(u_0^T) z'} (-\xi)
\end{split}
\end{equation*}  
we have
$$
\mathscr{L}^0 z^\sharp = \phi_1 + \phi_2 - \Big\{\mathscr{L}^0\big(\overline{\Pi}_{z'} X\big) + \Pi_{a (u_0^T) z'} (-\xi)\Big\},
$$
with initial condition
$$
z^\sharp (0) = z_0 - \overline{\Pi}_{z' (0)} X \in \mathcal{C}^{\alpha}.
$$
We can then apply Lemmas \ref{lem:a+b} and \ref{LemFixedPointSetting} to get the estimate \eqref{EqEstimateParac} from estimate \eqref{EqEstimate}. The estimate for $z^{\sharp} - y^{\sharp}$ is obtained by the same argument.
\end{Dem}

\bigskip
\bigskip

We now have in hands all we need to prove Theorem \ref{ThmMain} in its abstract form. Recall that the map 
$$
\Phi({\sf u}) = \Phi(u,u') =: {\sf v} =: (v,v')
$$ 
was defined at the end of Section \ref{SubsectionFixedPtSetting}, and note as a preliminary remark that $\big\|v^\sharp_{t=0}\big\|_{C^\alpha} \leq C(u_0)\,\big(1+\|X\|_{C^\alpha}\big)$, for a positive constant depending only on $\|u_0\|_{C^\alpha}$.

\medskip

\begin{thm}
\label{ThmAbstract}
The map $\Phi$ is a contraction from $\frak{B}_T(\lambda)$ into itself, for $\lambda$ large enough and $T$ small enough.
\end{thm}

\medskip

\begin{Dem}
Given the definition of $\Phi$ and the size estimates proved in Lemma \ref{LemVSharp}, Lemma \ref{LemmaV'} and Corollary \ref{Corollary}, we see that the bound
$$
\big\|\Phi({\sf u})\big\|_{\alpha,\beta} \lesssim C_1 + C_2\big(\|{\sf u}\|_{\alpha,\beta}\big)\,T^\frac{\alpha-\beta}{2}
$$
holds for some positive constants $C_i$ depending on the data, and $C_2$ depending also on $\|{\sf u}\|_{\alpha,\beta}$. The set $\frak{B}_T(\lambda)$ is then sent into itself by $\Phi$ for an adequate choice of parameters $\lambda$ and $T$.

\ssk

The map inherits its contracting character on $\frak{B}_T(\lambda)$ from the Lipschitz estimates given in Lemma \ref{LemVSharp}, Lemma \ref{LemmaV'} and Corollary \ref{Corollary}, and the fact that $v^\sharp_{t=0}$ is fixed, and depends only on $u_0$, for functions $v^\sharp$ built from $\Phi$.
\end{Dem}

\medskip

As the map $\Phi$ depends continuously on the enhanced noise $\widehat\xi := \big(\xi, \Pi(X,\xi)\big)\in C^{\alpha-2}\times C^{2\alpha-2}$, and the contracting character of $\Phi$ is locally uniform on $\widehat\xi$, its \textit{fixed point depends continuously on} $\widehat\xi$; we denote it by $\frak{I}\big(\widehat\xi\,\big)$. Given a zero spatial mean smooth function $\zeta$ on the $2$-dimensional torus, denote by $Z$ the solution to the equation $-\Delta Z = \zeta$. The function $\frak{I}$ extends the solution map $I$ that associates to a smooth 'noise' $\zeta$ the solution to the well-posed quasilinear equation 
$$
\partial_t v - a(v)\Delta v = g(v)\zeta, \quad v_{t=0} = u_0\in C^\alpha,
$$
in the sense that 
$$
I(\zeta) = \frak{I}\big(\zeta,\Pi(Z,\zeta)\big),
$$
for any smooth noise. (The fact that all these functions can be defined on the same time interval is part of the claim.)

\bigskip

\subsection{Renormalisation}
\label{SectionRenormalisation}

The above analysis of the singular quasilinear equation \eqref{EqQuasigPAM} requires that we start from the data of $\xi\in C^{\alpha-2}$ and $\Pi(X,\xi)\in C^{2\alpha-2}$. For a typical realization of space white noise, $X$ is only in $C^\alpha$ and one cannot make sense of the resonant term $\Pi(X,\xi)$ on a purely analytical basis. Gubinelli, Imkeller and Perkowski first showed in \cite{GIP} that if $\xi^\epsilon := P_\epsilon\xi$, and $X^\epsilon$ solves the equation $-\Delta X^\epsilon := \xi^\epsilon$, then there exists diverging constants $c^\epsilon$ such that 
$$
\Pi\big(X^\epsilon,\xi^\epsilon\big) - c^\epsilon
$$ 
converges almost surely in $C^{2\alpha-2}$ to some limit element, denoted by $\Pi(X,\xi)$. At the same time, identities \eqref{EqRenormalisation1} and \eqref{EqRenormalisation2} make it clear that 
$$
\frak{I}\Big(\xi^\epsilon,\Pi\big(X^\epsilon,\xi^\epsilon\big) - c^\epsilon\Big) = I^\epsilon\big(\xi^\epsilon\big),
$$
where $I^\epsilon$ stands for the solution map that associates to a smooth 'noise' $\zeta$ the solution to the well-posed quasilinear equation 
$$
\partial_t v - a(v)\Delta v = g(v)\zeta - c^\epsilon\left\{ \frac{g'g}{a}-a'\Big(\frac{g}{a}\Big)^2\right\}(v), \quad v_{t=0} = u_0.
$$
Theorem \ref{ThmMain} follows as a consequence of the continuity of the map $\frak{I}$ and the almost sure convergence of $\Big(\xi^\ep,\Pi\big(X^\epsilon,\xi^\epsilon\big) - c^\epsilon\Big)$ to $\big(\xi,\Pi(X,\xi)\big)$ in $C^{\alpha-2}\times C^{2\alpha-2}$, together with the remark that for $(u,u')\in {\sf C}_{\alpha,\beta}(X)$, we have the bound
$$
\|u\|_{\mcC^\alpha} \lesssim \|u'\|_{\mcC^\beta}\|X\|_{C^\alpha} + \big\|u^\sharp\big\|_{\mcC^\alpha},
$$
a consequence of the elementary Lemma \ref{lem:l4}.

\bigskip

\appendix

\section{Elementary side results}
\label{Appendix}

We collect in this Appendix a number of elementary side results, together with their proofs, to make this work self-contained. They can all be found somewhere else in some form or another. As a warm-up, we start with some elementary claims, used in the main body of the text. Recall the classical notation $\Delta_i$ for the Fourier multipliers used in Littlewood-Paley decomposition; refer to Bahouri, Chemin and Danchin's textbook \cite{BCD} for the basics on this subject.

\medskip

\begin{lem}  
\label{lem:aux}
\begin{enumerate}
    \item[{\sf 1.}] For $u$ in the spatial H\"older space $C^\alpha$, we have
    $$
	\big\| g (u) \big\|_{C^\alpha} \leq \| g \|_{C^1} \big(1 + \| u\|_{C^\alpha}\big),    
    $$
	and 
	$$
	\big\| g (u) - g (v) \big\|_{C^\alpha} \leq \| g \|_{C^2} \big(1 + \| u \|_{\alpha}\big) \| u - v \|_{C^\alpha}.
	$$   \vspace{0.1cm}
    
    \item[{\sf 2.}] For $u$ in the parabolic H\"older space $\mcC^\alpha$, and $0< \beta\leq \alpha$, we have
 	$$
 	\big\| g (u) \big\|_{\mcC^\beta} \lesssim T^{\frac{\alpha - \beta}{2}} \| g \|_{C^1} \big(1 + \| u \|_{\mcC^\alpha}\big) + \big\| g (u_0)\big\|_{C^\beta},
 	$$   
	and if $u_{t=0}=v_{t=0}$, then 
	$$
	\big\| g (u) - g (v)\big\|_{\mcC^\beta} \lesssim T^{\frac{\alpha - \beta}{2}} \| g\|_{C^2}  \big(1 + \| u \|_{\mcC^\alpha}\big) \| u - v\|_{\mcC^\alpha}.
	$$
\end{enumerate}
\end{lem}

\medskip

\begin{Dem}
The Lipschitz estimate of point {\sf 1} comes by writing 
$$
\big(g(u)-g(v)\big)(x) - \big(g(u)-g(v)\big)(y)
$$
as the boundary term of the integral of the derivative of the function $s \in [0,1] \mapsto\big(g(u)-g(v)\big)\big(y+s(x-y)\big)$.

\ssk

To see point {\sf 2}, take any function $h\in\mcC^\beta$ and start from the following two estimates
$$
\big\| \Delta_i \big\{h(t) - h(0)\big\} \big\|_{L^{\infty}} \lesssim \left\{ \begin{array}{l}
       2^{- i \alpha} \, \|h\|_{C_T \mathcal{C}^{\alpha}}   \\
       T^{\frac{\alpha}{2}} \, \|h\|_{C^{\alpha / 2}_T L^{\infty}}
     \end{array} \right.
$$
to get by interpolation of the two upper bounds the estimate
\begin{equation*}
\|h\|_{C_T C^\beta} \lesssim T^{\frac{\alpha - \beta}{2}} \|h\|_{C_T C^\alpha}^\varepsilon\, \|h\|^{1 - \varepsilon}_{C_T^\frac{\alpha}{2} L^{\infty}} + \big\|h(0)\big\|_{C^\beta}.
\end{equation*}
It follows that if $u_{t=0}=v_{t=0}$, then we have
\begin{equation*}
\begin{split}
\big\| g (u) - g (v) \big\|_{C_T \mathcal{C}^{\beta}} &\lesssim T^{\frac{\alpha -\beta}{2}} \, \big\| g (u) - g (v) \big\|_{C_T C^\alpha}^{\varepsilon} \big\| g (u) - g (v) \big\|^{1 - \varepsilon}_{C_T^{{\alpha/2}} L^{\infty}}   \\
&\lesssim T^{\frac{\alpha - \beta}{2}} \| g \|_{C^2}  \Big\{(1 + \| u \|_{C_T \mathcal{C}^{\alpha}}) \| u - v \|_{C_T \mathcal{C}^{\alpha}}   \\
&\hspace{2.5cm}+ \Big(1 + \| u\|_{C_T^{\alpha/2} L^{\infty}}\Big) \| u - v \|_{C_T^{\alpha/2} L^{\infty}}\Big\}.
\end{split}
\end{equation*}
The Lipschitz estimate of point {\sf 2} then comes as a consequence of the inequality
\begin{equation*}
\begin{split}
\big\| g (u) - g (v) \big\|_{C^{\beta / 2}_T L^{\infty}} &\leq T^{\frac{\alpha - \beta}{2}} \, \big\| g (u) - g (v) \big\|_{C^{\alpha / 2}_T  L^{\infty}}   \\
&\leq T^{\frac{\alpha - \beta}{2}} \| g \|_{C^2} \Big(1 + \|u \|_{C^{\alpha / 2}_T L^{\infty}}\Big) \| u - v \|_{C^{\alpha / 2}_T L^{\infty}}.
\end{split}
\end{equation*}
\end{Dem}

\medskip

The next result is a variation on Gubinelli, Imkeller and Perkowski's fundamental 'commutator' lemma; it is the first part of Lemma \ref{LemFixedPointSetting}. It also happens to be special case of a more general result, proved in Theorem 2 of \cite{BailleulBernicot2}. Recall that we work with $\alpha>\frac{2}{3}$.

\medskip

\begin{lem}
Let $f, g \in C^\beta$, and $a \in C^\alpha$ and $b \in C^{\alpha - 2}$ be given, such that $\Pi(a,b)$ is a well-defined element of $C^{2\alpha-2}$. Then
$$
\Big\| \Pi\big(\Pi_f a, \Pi_g b\big) - (f g) \Pi (a,b) \Big\|_{C^{2 \alpha + \beta - 2}} \lesssim \| f \|_{C^\beta} \| g \|_{C^\beta} \|a\|_{C^\alpha} \|b\|_{C^{\alpha - 2}}.
$$
\end{lem}

\medskip

\begin{Dem}
Denote by $K_{k, x} (z) := K_k (x - z)$ the convolution kernel of the Littlewood-Paley projector $\Delta_k$. We have
\begin{equation*}
\begin{split}
\Delta_k \Big(\Pi\big(\Pi_f a, \Pi_g b\big)\Big) (x) &= \int K_{k, x} (y) \Big(\Pi\big(\Pi_{f (y)}a, \Pi_{g (y)}b\big)\Big)(y)\,dy   \\
&\quad+ \int K_{k, x} (y) \Big(\Pi\big(\Pi_{f - f (y)}a, \Pi_{g (y)}b\big)\Big) (y)\,dy   \\
&\quad+ \int K_{k, x} (y) \Big(\Pi \big(\Pi_f a, \Pi_{g - g (y)}b\big)\Big) (y)\,dy   \\
&=: {\sf I}_1 + {\sf I}_2 + {\sf I}_3.
\end{split}
\end{equation*}
The first term gives
\begin{equation*}
\begin{split}
\Big(\Pi\big(\Pi_{f (y)}a, \Pi_{g (y)}b\big)\Big) (y) &= \sum_{i \sim j} \Delta_i \big(\Pi_{f (y)}a\big)(y) \Delta_j\big(\Pi_{g (y)}b\big)(y)   \\     
&= \sum_{i \sim j} f (y) g (y) \big(\Delta_i a\big) (y) \big(\Delta_j b\big) (y) = f (y) g (y) \, \Pi (a, b) (y)
\end{split}
\end{equation*}
so it remains to estimate ${\sf I}_2$ and ${\sf I}_3$ in the spatial H\"older space $C^{2\alpha + \beta - 2}$. We have
\begin{equation*}
\begin{split}
| {\sf I}_2 | &= \left| \int K_{k, x} (y) \sum_{i \sim j \gtrsim k} \Delta_i\big(\Pi_{f - f (y)}a\big) (y) g (y) \big(\Delta_j b\big) (y) \,dy\right|   \\
&\lesssim \| f \|_{C^\alpha} \|a\|_{C^\alpha} \|g\|_{L^{\infty}} \| b\|_{C^{\alpha - 2}} \int \big|K_{k, x}(y)\big| \left(\sum_{i \gtrsim k} 2^{- (2 \alpha+ \beta - 2) i}\right)\,dy;
\end{split}
\end{equation*}     
since $2 \alpha + \beta - 2 > 0$ we obtain in the end
$$
| {\sf I}_2 | \lesssim 2^{- (2 \alpha + \beta - 2) k} \| f \|_{C^\alpha} \| a\|_{C^\alpha} \| g \|_{L^{\infty}} \| b \|_{C^{\alpha - 2}}
$$
Similarly, we have 
$$
| {\sf I}_3 | \lesssim 2^{- (2 \alpha + \beta - 2) k} \| f \|_{L^{\infty}} \| a\|_{C^\alpha} \| g \|_{C^\beta} \|b\|_{C^{\alpha - 2}},
$$
which completes the proof.
\end{Dem}

\medskip

\begin{lem}
\label{lem:l3}
Given $f \in C^{2 \beta}$, $g \in C^\beta$, and $h \in C^{\alpha - 2}$, with regularity exponents in $(0,2)$, we have
$$
\big\| f \, \Pi_g h - \Pi_{f g} h \big\|_{C^{\beta + \alpha - 2}} \lesssim \|f\|_{C^{2 \beta}} \| g \|_{C^\beta} \| h \|_{C^{\alpha - 2}}.
$$
\end{lem}

\medskip

\begin{Dem}
We have
\begin{equation*}
\begin{split}
f \, \Pi_g h &= \Pi_f (\Pi_g h) + \Pi_{\Pi_g h} (f) + \Pi\big(f, \Pi_gh\big)   \\
				  &= \big(\Pi_f (\Pi_g h) - \Pi_{f g}h\big) + \Pi_{f g}h + \Pi_{\Pi_g h}(f) + \Pi \big(f, \Pi_g h\big),
\end{split}
\end{equation*}
where
$$
\big\| \Pi_f \big(\Pi_g h\big) - \Pi_{f g}h \big\|_{C^{\beta + \alpha - 2}} \lesssim \|f \|_{L^{\infty}}\, \| g \|_{C^\beta} \, \| h \|_{{C^{\alpha - 2}}},
$$
due to Proposition 23 of \cite{BailleulBernicot2}, and
\begin{equation*}
\begin{split}
&\big\| \Pi \big(f, \Pi_g h\big)\big\|_{C^{2\beta + \alpha - 2}} \lesssim \| f \|_{C^{2\beta}} \,\| g \|_{L^{\infty}}\, \| h \|_{C^{\alpha - 2}}   \\
&\big\| \Pi_{\Pi_g h} f \big\|_{C^{2\alpha - 2}} \lesssim \| f \|_{C^\alpha} \, \| g\|_{L^{\infty}} \,\| h \|_{C^{\alpha - 2}}.
\end{split}
\end{equation*}
\end{Dem}

\medskip

The next proposition gives the second part of Lemma \ref{LemFixedPointSetting}; recall $\mathscr{L}^0 = \partial_t - a(u_0)\Delta$.

\medskip

\begin{prop}
\label{prop:2}
Given $u'\in\mcC^\beta$, we have the following continuity result for a commutator
$$
\Big\|\mathscr{L}^0 \big(\overline\Pi_{u'} X\big) - \Pi_{a (u^T_0) u'} (-\Delta X) \Big\|_{C_T\mathcal{C}^{\alpha + \beta - 2}} \lesssim \big(1 + T^{- \frac{2\beta-\alpha}{2}} \| u_0\|_{C^\alpha}\big) \| u' \|_{\mcC^\beta} \| X \|_{C^\alpha}.
$$
\end{prop}

\medskip

\begin{Dem}
Recall first from Lemma 5.1 in  \cite{GIP} that
\begin{equation}
\label{eq:23}
\big\| \overline\Pi_{u'} (\Delta X) - \Pi_{u'} (\Delta X) \big\|_{C_T C^{\alpha + \beta - 2}} \lesssim \| u' \|_{C^{\beta / 2}_TL^{\infty}} \| X \|_{C^\alpha}.
\end{equation}
Then we write
\begin{equation*}
\begin{split}
\mathscr{L}^0 \big(\overline\Pi_{u'}X\big) - \Pi_{a (u^T_0) u'} (-\Delta X) &= \Big\{\mathscr{L}^0 \big(\overline\Pi_{u'}X\big) + a (u^T_{^{} 0}) \overline\Pi_{u'} (\Delta X)\Big\}   \\
&\quad+ a (u^T_0) \Big\{\Pi_{u'} (\Delta X) - \overline\Pi_{u'} (\Delta X)\Big\}   \\
&\quad+ \Big\{\Pi_{a (u^T_0) u'} (\Delta X) - a (u^T_0) \Pi_{u'} (\Delta X)\Big\},
\end{split}
\end{equation*}
and observe that the second term on the right hand side can be estimated with inequality \eqref{eq:23} 
$$
\Big\| a (u^T_0) \Big\{\Pi_{u'} (\Delta X) - \overline\Pi_{u'} (\Delta X)\Big\}\Big\|_{C^{\alpha + \beta - 2}} \lesssim \big(1 + \| u_0 \|_{\alpha}\big) \| u'\|_{C^{\beta / 2}_T L^{\infty}} \| X \|_{C^\alpha},
$$
and that the third term can be taken care of by Lemma \ref{lem:l3}
$$
\big\| \Pi_{a (u^T_0) u'} (\Delta X) - a (u^T_0) \Pi_{u'} (\Delta X) \big\|_{C^{\alpha + \beta - 2}} \lesssim \big(1 + T^{- \frac{2\beta-\alpha}{2}} \| u_0\|_{C^\alpha}\big) \| u' \|_{\mcC^\beta} \| X \|_{C^\alpha}.
$$
We now estimate the first term. Since $X$ does not depend on $t$, 
$$
\partial_t \big(\overline\Pi_{u'} X\big) = \sum_i \partial_t (S_{i - 1} Q_i u') \Delta_i X.
$$
The spatial Fourier transform of $\partial_t \big(S_{i - 1} Q_i u'\big) \Delta_i X$ is localized in an annulus of size $2^i$, we obtain from estimate (32) in \cite{GIP} that
\begin{equation*}
\begin{split}
\big\| \partial_t (S_{i - 1} Q_i u')\big\|_{C_T L^{\infty}} &= \big\|\partial_t (Q_i S_{i - 1} u')\big\|_{C_T L^{\infty}}   \\
&\lesssim 2^{- (\beta - 2) i} \big\| S_{i - 1} u' \big\|_{C^{\beta / 2}_T L^{\infty}}   \\
&\lesssim 2^{- (\beta - 2) i} \| u' \|_{C^{\beta / 2}_T L^{\infty}},
\end{split}
\end{equation*}
so
\begin{equation}    
\label{eq:21}
\big\| \partial_t (\overline\Pi_{u'}X\big)\big\|_{C_T C^{\alpha + \beta - 2}} \lesssim \| u' \|_{\mcC^\beta} \| X \|_{C^\alpha}.
\end{equation}
We have, on the other hand, 
$$
\Delta \big(\overline\Pi_{u'}X\big) - \overline\Pi_{u'} (\Delta X) = \overline\Pi_{\Delta u}X - 2 \overline\Pi_{\nabla u'} (\nabla X),
$$
with
$$
\big\| Q_i S_{i - 1} \Delta u' \big\|_{C_T L^{\infty}} \leq \big\| S_{i - 1} \Delta u' \big\|_{C_T L^{\infty}} \lesssim 2^{- (\beta - 2) i} \| u' \|_{C_T C^\beta}
$$
and  
$$
\big\| Q_i S_{i - 1} \nabla u' \big\|_{C_T L^{\infty}} \lesssim 2^{- (\beta - 1) i} \| u' \|_{C_T C^\beta}.
$$
Altogether, this gives
\begin{equation}
\label{eq:22}
\big\| \Delta (\overline\Pi_{u'}X) - \overline\Pi_{u'} (\Delta X) \big\|_{C_T C^{\alpha + \beta - 2}} = \| \overline\Pi_{\Delta u} (X) - 2 \overline\Pi_{\nabla u'} (\nabla X) \|_{C_T C^{\alpha + \beta - 2}} \lesssim \| u' \|_{\mcC^\beta} \| X \|_{C^\alpha},
\end{equation}
so we deduce from \eqref{eq:21} and \eqref{eq:22} that
\begin{equation*}
\begin{split}
\Big\| \mathscr{L}^0 (\overline\Pi_{u'} (X)) &- a (u^T_0) \overline\Pi_{u'} (-\Delta X) \Big\|_{C_T C^{\alpha + \beta - 2}}   \\
&= \Big\| \partial_t \big(\overline\Pi_{u'}X\big) - a (u^T_0) \Big\{\Delta(\overline\Pi_{u'} (X)) - \overline\Pi_{u'} (\Delta X)\Big\} \Big\|_{C_T C^{\alpha + \beta - 2}}   \\
&\lesssim \big(1 + \| a (u^T_0) \|_{C^\beta}\big) \| u' \|_{\mcC^\beta} \| X \|_{C^\alpha},
\end{split}
\end{equation*}
which concludes the proof.
\end{Dem}

\medskip

\begin{lem}
\label{lem:l4}
Given $f \in \mcC^\beta$ and $g \in C^\alpha$, we have
$$
\big\| \overline\Pi_fg \big\|_{\mcC^\alpha} \lesssim \| f \|_{\mcC^\beta} \| g \|_{C^\alpha}.
$$
\end{lem}

\medskip

\begin{Dem}
Let work here with the canonical heat operator $\mathscr{L} := \partial_t - \Delta$. We have from the classical Schauder estimates
$$
\big\|\overline\Pi_fg\big\|_{\mcC^\alpha} \lesssim \big\|\overline\Pi_{f(0)}g\big\|_{C^\alpha} + \big\| \mathscr{L}\big(\overline\Pi_fg\big)\big\|_{C_T C^{\alpha - 2}},
$$
and the rough bound $\big\|\overline\Pi_{f(0)}g\big\|_{C^\alpha} \leq \|f\|_{C_T L^{\infty}} \|g\|_{C^\alpha}$. Next, we write, with $\mathscr{L}g = -\Delta g$,
$$
\mathscr{L} \big(\overline\Pi_fg\big) = \Big\{\mathscr{L}\big(\overline\Pi_fg\big) - \overline\Pi_f (-\Delta g)\Big\} - \overline\Pi_f (\Delta g),
$$
and use commutator Lemma 5.1 of \cite{GIP} to get
$$
\Big\| \mathscr{L}\big(\overline\Pi_fg\big) - \overline\Pi_f (-\Delta g) \Big\|_{C_T\mathcal{C}^{\beta + \alpha - 2}} \lesssim \| f \|_{\mcC^\beta}\| g \|_{C^\alpha},
$$
and
$$
\Big\| \overline\Pi_f (-\Delta g)\Big\|_{C_T C^{\alpha - 2}} \lesssim \|f \|_{C_T L^{\infty}} \| g \|_{C^\alpha}.
$$
\end{Dem}

\bigskip
\bigskip

\noindent \textcolor{gray}{$\bullet$} {\sf I. Bailleul} - {\small IRMAR, Universit\'e de Rennes 1, France.}   \vspace{0.1cm}

\vspace{-0.1cm}\noindent \hspace{0.2cm} {\it ismael.bailleul@univ-rennes1.fr}   \vspace{0.3cm}

\noindent \textcolor{gray}{$\bullet$} {\sf A. Debussche} - {\small IRMAR, \'Ecole Normale Sup\'erieure de Rennes, France.}   \vspace{0.1cm}

\vspace{-0.1cm}\noindent \hspace{0.2cm} {\it arnaud.debussche@ens-rennes.fr}   \vspace{0.3cm}

\noindent \textcolor{gray}{$\bullet$} {\sf M. Hofmanov\'a} - {\small Institute of Mathematics, Technical University Berlin, Germany.}   \vspace{0.1cm}

\vspace{-0.1cm}\noindent \hspace{0.2cm} {\it hofmanov@math.tu-berlin.de}

\end{document}